# Calculating the interaction index of a fuzzy measure: A polynomial approach based on sampling

Inmaculada Gutiérrez [a,*], Javier Castro [a,b], Daniel Gómez [a,b], Rosa Espínola [a,b]

[a] *Faculty of Statistics, Complutense University, Madrid, Spain*
[b] *Instituto de Evaluación Sanitaria, Complutense University, Madrid, Spain*

## Abstract

In this paper we address the problem of fuzzy measures index calculation. On the basis of fuzzy sets, Murofushi and Soneda proposed an interaction index to deal with the relations between two individuals. This index was later extended in a common framework by Grabisch. Both indices are fundamental in the literature of fuzzy measures. Nevertheless, the corresponding calculation still presents a highly complex problem for which no approximation solution has been proposed yet. Then, using a representation of the Shapley based on orders, here we suggest an alternative calculation of the interaction index, both for the simple case of pairs of individuals, and for the more complex situation in which any set could be considered. This alternative representation facilitates the handling of these indices. Moreover, we draw on this representation to define two polynomial methods based on sampling to estimate the interaction index, as well as a method to approximate the generalized version of it. We provide some computational results to test the goodness of the proposed algorithms.
© 2022 The Author(s). Published by Elsevier B.V. This is an open access article under the CC BY license (http://creativecommons.org/licenses/by/4.0/).

## 1. Introduction

In the second half of the 20th century, the measure theory [1] was increasingly popular. This field has been crucial in the development of several scientific areas, as the fuzzy measures introduced by Sugeno [2], the necessity and possibility measures introduced by Zadeh [3,4], or the capacity measures introduced by Choquet [5], among many others.

We focus on fuzzy measures. These monotonic set functions define a type of non-additive measures encompassing many measures such as belief, possibility, necessity or plausibility measures [6–8]. Fuzzy measures have been widely analyzed [33,9] and applied in many fields [10]. For example, in the approximate reasoning context are useful to make decisions, choose proper methods and logical operators for implications and connectives [11]. These are also essential in the context of decision theory [12], multi-criteria decision aid [13], when dealing with the Choquet integral [14] or to face aggregation problems [11,15]. One of the strengths of fuzzy measures and most applied features is their ability to model interactions between elements. Even, in the context of multi-criteria problem, those functions can be used to represent interactions among criteria [16]. Unfortunately, the practical implementation of fuzzy measures is quite difficult: $2^n$ coefficients are needed to define a fuzzy measure over $n$ individuals. Also the related semantic involved is difficult to understand and interpret [17].

---

\* Corresponding author.
*E-mail addresses:* inmaguti@ucm.es (I. Gutiérrez), jcastroc@estad.ucm.es (J. Castro), dagomez@estad.ucm.es (D. Gómez), rosaev@estad.ucm.es (R. Espínola).

Many researchers have focused on the characterization, interpretation and representation of fuzzy measures. The question arises whether it is really worth identifying complete fuzzy measures, as in real problems, the amount of interactions is limited. In these terms, it is worth mentioning the work of Beliakov and Wu [18], who define the $k$-interactivity, or the proposal of Grabisch [19] on the consideration of $k$-order additive fuzzy measures, in addition to other techniques to learn about fuzzy measures [9,20]. Defining multiple indices encompasses multiple approaches, in an attempt to reduce the complexity. In 1953, Shapley defined an importance index: the Shapley value [21]. It is essential in Game Theory, and has been adapted to the field of fuzzy measures, where it can be seen as the "importance" of each singleton. The Shapley value has been deeply studied from a theoretical perspective [22,23]; however, there is an issue on the complexity: its computation is a $NP$-problem. Thus, despite its importance, from the best of our knowledge, until not that long ago, there were only a few approaches to approximate it. Furthermore, these references were focused on some specific problems, as when the population has specific properties [24,25]. We found more works on the approximation of the Shapely value, as [26] or [27]. In 2009, Castro et al. proposed a method based on sampling to approximate the Shapley value [28]. It was generalized in [29] by a process of stratified random sampling with optimum allocation. That method is at least as good as the first one. However, random sampling based method should not be underestimated; it is clear and simple, with easy understanding and quick computation, which allows a first sight of the refinement with stratification.

Besides the importance of each singleton represented by the Shapley value, it could be interesting to analyze the interactions between elements [30–32]. The problem is that, when using fuzzy measures to model interactions among elements, the number of coefficients exponentially grows with the number of individuals, likewise the great computational complexity of assigning a value to all the coalitions. Grabisch proposed a technique to estimate the importance of features and their interaction [33]. On the basis of fuzzy sets [34,35], Murofushi and Soneda characterized the interactions among two individuals [36]: given a coalition, the interaction index quantifies its average contribution, considering all the subsets it is part of. Besides its powerful at modeling relations between variables, items or criteria, the interaction index is an essential tool at representation. Grabisch generalized this index to deal with bigger sets [19]. As with the Shapley value, the approximation of the real value of these indices has not been studied in depth for general scenarios, due to its great complexity.

Then, we propose an alternative calculation of this interaction index, for both the simple case related to pairs of elements and the extension to deal with any set. This new formulation of the interaction index (and in general of the representation index), provides us a different and intuitive understanding of this index, based on orders. Regarding the idea of sampling introduced in [28,29], we propose two methods to approach the value of the interaction index. Both are based on sampling, where the second one is a refinement of the very first, but also using stratified random sampling with optimum allocation. In our view, it is important to maintain both proposals, not only the refinement, as the first one, based on random sampling, is easier to compute, and facilitates the comprehension and understanding of the stratified version. As mentioned before, the second approach will perform better in some easily recognizable situations. We will show that both methods have desirable statistical properties. In our humble opinion, being able to do this calculation motivates the usefulness and justifies the importance of this new characterization.

The remainder of the paper is organized as follows. In Section 2, we provide a detailed explanation about the interaction index and an alternative calculation of it. In Section 3 we propose an algorithm based on sampling theory for the estimation of the interaction index. A refinement of that algorithm, based on stratification and maximum allocation, is introduced in Section 4. In Section 5 we evaluate the performance of the proposed algorithms. We first calculate, in polynomial time, the real value of the interaction index in some basic examples. Then, we compare the obtained real values with the estimations provided by both the *ApproInteraction* and *Stratified ApproInteraction* algorithms. We carry on with some conclusions and final remarks in Section 6.

## 2. A new representation of the interaction index based on orders

On the following we will refer to the set $X = \{1, \ldots, n\}$, whose individuals may refer to the attributes or criteria of a decision problem [37,38] or the players of a cooperative game [39], for example. Specifically, we work with fuzzy measures [2,19], monotonic set functions as $\mu : 2^X \to [0, 1]$ with $\mu(\emptyset) = 0$. We work with the Shapley value [21], calculated for the individual $i \in X$ as $Sh_i(\mu) = \sum_{k=0}^{n-1} \gamma_k \sum_{\substack{K \subset X \setminus \{i\} \\ |K|=k}} (\mu_{iK} - \mu_K)$, being $\gamma_k = \frac{(n-k-1)!k!}{n!} = \frac{1}{\binom{n-1}{k}n}$ (where $\mu_i, \mu_{ij}, \ldots, \mu_K, \mu_{iK}$ denote $\mu(\{i\})$, $\mu(\{i, j\})$, $\ldots$, $\mu(K)$, $\mu(\{i\} \cup K)$, respectively). We also consider an alternative calculation of the Shapley value based on orders: being $\pi(X)$ the set of all possible permutations of the elements of $X$; $o \in \pi(X)$ one of these orders; and $Pred_i(o)$ the set of predecessors of $i \in X$ in the order $o$, the Shapley value of $i$ can be rewritten as $Sh_i(\mu) = \frac{1}{n!} \sum_{o \in \pi(X)} [\mu(Pred_i(o) + \{i\}) - \mu(Pred_i(o))]$.

Beyond the Shapley value Murofushi and Soneda [36] proposed an interaction index to model the interactions between pairs of elements. Based on concepts of multi-attribute utility theory, it is a powerful representation tool characterized as $I_{ij}(\mu) = \sum_{k=0}^{n-2} \zeta_k \sum_{\substack{K \subset X \setminus \{i,j\} \\ |K|=k}} (\mu_{ijK} - \mu_{iK} - \mu_{jK} + \mu_K)$, where $\zeta_k = \frac{(n-k-2)!k!}{(n-1)!} = \frac{1}{\binom{n-2}{k}(n-1)}$. A deep analysis of $I_{ij}$ can be found in [40]. As we will demonstrate, analogously to the Shapely value, $I_{ij}(\mu)$ can be seen as an average of the added value obtained by considering the individuals $i$ and $j$ together in the same coalition, depending on $\mu$. To do so, we will suggest an alternative calculation of the $I_{ij}$ based on orders.

To deal with general scenarios, Grabisch defined the representation index (or interactions representation) to model interactions among individuals of any subset $T \subset X$ [19]: $I_T(\mu) = \sum_{k=0}^{n-|T|} \zeta_k^{|T|} \sum_{\substack{K \subset X \setminus T \\ |K|=k}} \sum_{L \subset T} (-1)^{|T|-|L|} \mu_{LK}$, where $\zeta_k^{|T|} = \frac{(n-k-|T|)!k!}{(n-|T|+1)!} = \frac{1}{\binom{n-|T|}{k}(n-|T|+1)}$. It is known as representation index or interactions representation.

Note that, for $T = \{i\}$, $I_T$ recovers the Shapley index $Sh_i(\mu)$. Similarly, it recovers the interaction index $I_{ij}$ for $T = \{i, j\}$. So, when considering sets with more than two elements, the characterization of $I_T$ is consistent with the one devoted to interactions between two individuals, $I_{ij}$.

Inspired by the characterization of the Shapley value based on orders, we propose an alternative calculation of the interaction index $I_{ij}$ (and also of $I_T$) based on orders, rewriting $I_{ij}$ in terms of all possible orders of the individuals. This characterization provides an overview of $I_{ij}$ (or $I_T$) as an average of the added value obtained by considering two individuals (or the individuals of $T$) together in the same coalition, depending on the involved fuzzy measure.

**Theorem 1.** *The interaction index $I_{ij}$ of the fuzzy measure $\mu : 2^X \to [0, 1]$ can be calculated as*

$$I_{ij}(\mu) = \sum_{o \in \pi(X \setminus \{i\})} \frac{1}{(n-1)!} [\mu(Pred_j(o) \cup \{i\} \cup \{j\}) \tag{1}$$
$$- \mu(Pred_j(o) \cup \{i\}) - \mu(Pred_j(o) \cup \{j\}) + \mu(Pred_j(o))]$$

**Proof.** Let

$$I_{ij}(\mu) = \sum_{k=0}^{n-2} \frac{(n-k-2)!k!}{(n-1)!} \sum_{\substack{K \subset X \setminus \{i,j\} \\ |K|=k}} (\mu_{ijK} - \mu_{iK} - \mu_{jK} + \mu_K)$$

$$= \sum_{k=0}^{n-2} \frac{(n-k-2)!k!}{(n-1)!} \sum_{\substack{K \subset X \setminus \{i,j\} \\ |K|=k}} (\mu_{ijK} - \mu_{iK}) - \sum_{k=0}^{n-2} \frac{(n-k-2)!k!}{(n-1)!} \sum_{\substack{K \subset X \setminus \{i,j\} \\ |K|=k}} (\mu_{jK} - \mu_K)$$

1. We define $\mu' : 2^{X'} \to [0, 1]$, such that $\forall S \subseteq X'$, $\mu'(S) = \mu(S \cup \{i\})$, being $X' = X \setminus \{i\}$.

  Then, the element $\sum_{k=0}^{n-2} \frac{(n-k-2)!k!}{(n-1)!} \sum_{\substack{K \subset X \setminus \{i,j\} \\ |K|=k}} (\mu_{ijK} - \mu_{iK})$ is $Sh_j(\mu')$, i.e. the capacity of $j$ when the individ-

  ual $i$ is in all its possible coalitions.

2. Similarly, we define $\mu'' : 2^{X''} \longrightarrow [0,1]$, such that $\forall S \subseteq X''$, $\mu''(S) = \mu(S)$, being $X'' = X' = X \setminus \{i\}$.

 Then, the element $\sum_{k=0}^{n-2} \frac{(n-k-2)!k!}{(n-1)!} \sum_{K \subset X \setminus \{i,j\}} (\mu_{jK} - \mu_K)$ is the Shapley index of the individual $j$ on the function $\mu''$. This is, the capacity of $j$ when it is never in a coalition with the individual $i$.

Because of points (1) and (2), we can rewrite the interaction index as

$$I_{ij}(\mu) = Sh_j(\mu') - Sh_j(\mu'') \tag{2}$$

where $Sh_j(\mu')$ (or $Sh_j(\mu'')$) is the Shapley index of the individual $j$ on $\mu'$ (or $\mu''$), in a coalition with all the elements of $X'$. Remember $X' \equiv X''$.

Due to the characterization of the Shapley index based on orders, we can rewrite the equation (2) as

$$
\begin{aligned}
I_{ij}(\mu) &= \sum_{o \in \pi(X')} \frac{1}{(n-1)!} \Big[ \mu'\big(Pred_j(o) \cup \{j\}\big) - \mu'\big(Pred_j(o)\big) \Big] \\
&\quad - \sum_{o \in \pi(X')} \frac{1}{(n-1)!} \Big[ \mu''\big(Pred_j(o) \cup \{j\}\big) - \mu''\big(Pred_j(o)\big) \Big] \\
&= \sum_{o \in \pi(X \setminus \{i\})} \frac{1}{(n-1)!} \Big[ \mu\big(Pred_j(o) \cup \{i\} \cup \{j\}\big) - \mu\big(Pred_j(o) \cup \{i\}\big) \Big] \\
&\quad - \sum_{o \in \pi(X \setminus \{i\})} \frac{1}{(n-1)!} \Big[ \mu\big(Pred_j(o) \cup \{j\}\big) - \mu\big(Pred_j(o)\big) \Big] \\
&= \sum_{o \in \pi(X \setminus \{i\})} \frac{1}{(n-1)!} \Big[ \mu\big(Pred_j(o) \cup \{i\} \cup \{j\}\big) \\
&\quad - \mu\big(Pred_j(o) \cup \{i\}\big) - \mu\big(Pred_j(o) \cup \{j\}\big) + \mu\big(Pred_j(o)\big) \Big]
\end{aligned}
$$

**Theorem 2.** *The representation index $I_T$ of the fuzzy measure $\mu : 2^X \to [0,1]$ can be calculated as*

$$I_T(\mu) = \sum_{o \in \pi(X \setminus \{t_2, \ldots, t_{|T|}\})} \frac{1}{(n - |T| + 1)!} [\chi(o)_T] \tag{3}$$

*where* $\quad \chi(o)_T = \sum_{L \subset T} (-1)^{|T| - |L|} \mu\big(Pred_{t_1}(o) \cup L\big)$

**Proof.** Being $L \in \{t_2, \ldots, t_{|T|}\}$, we define $\mu^L : 2^{X'} \longrightarrow [0,1]$, such that $\forall S \subseteq X'$, $\mu^L(S) = \mu(S \cup L)$, being $X' = X \setminus \{t_2, \ldots, t_{|T|}\}$. Then,

$$Sh_{t_1}(\mu^L) = \sum_{k=0}^{n-|T|} \zeta_k^T \sum_{\substack{K \subset X' \setminus \{t_2, \ldots, t_{|T|}\} \\ |K| = k}} \left( \mu_{t_1 K}^L - \mu_K^L \right) \tag{4}$$

*where* $\quad \zeta_k^T = \dfrac{(n - k - |T|)! k!}{(n - |T| + 1)!} = \dfrac{1}{\binom{n - |T|}{k}(n - |T| + 1)}$

Then,

$$
\begin{aligned}
I_T(\mu) &= \sum_{k=0}^{n-|T|} \zeta_k^{|T|} \sum_{\substack{K \subset X \setminus T \\ |K|=k}} \sum_{L \subset T} (-1)^{|T|-|L|} \mu_{LK} \\
&= \sum_{k=0}^{n-|T|} \zeta_k^{|T|} \sum_{\substack{K \subset X' \setminus \{t_2,\ldots,t_{|T|}\} \\ |K|=k}} \sum_{L \subset T \setminus t_1} (-1)^{|T|-|L|-1} (\mu_{LKt_1} - \mu_{LK}) \\
&= \sum_{L \subset T \setminus t_1} (-1)^{|T|-|L|-1} \sum_{k=0}^{n-|T|} \zeta_k^{|T|} \sum_{\substack{K \subset X' \setminus \{t_2,\ldots,t_{|T|}\} \\ |K|=k}} (\mu_{Kt_1}^L - \mu_K^L) \\
&= \sum_{L \subset T \setminus t_1} (-1)^{|T|-|L|-1} Sh_{t_1}(\mu^L)
\end{aligned}
$$

(5)

Due to the characterization of the Shapley index based on orders and the definition of $Sh_{t_1}(\mu^L)$ given in equation (4) we can rewrite equation (5) as:

$$
\begin{aligned}
I_T(\mu) &= \sum_{L \subset T \setminus t_1} (-1)^{|T|-|L|-1} \sum_{o \in \pi(X')} \frac{1}{(n-|T|+1)!} (\mu^L(Pred_{t_1}(o) \cup \{t_1\}) - \mu^L(Pred_{t_1}(o))) \\
&= \sum_{o \in \pi(X \setminus \{t_2,\ldots,t_{|T|}\})} \frac{1}{(n-|T|+1)!} \sum_{L \subset T} (-1)^{|T|-|L|} \mu(Pred_{t_1}(o) \cup L)
\end{aligned}
$$

## 3. On the calculation of the estimation of the interaction index

As mentioned, $I_{ij}$ is an essential tool in the field of fuzzy measures whose calculation is usually not easy. In fact, $I_{ij}$ can be calculated in polynomial time only for a few simple fuzzy measures. Also, to the best of our knowledge, there is not much research on feasible approximations of this index.

Then we propose a polynomial-time algorithm to estimate the value of $I_{ij}$. We also provide some desirable properties of this estimated value. To do so, we define a unique sampling process for every pair $i, j \in X$, named *ApproInteraction* Algorithm, is defined below (see Algorithm 1). To simplify the notation we denote the estimated value of $I_{ij}(\mu)$ as $\hat{I}_{ij} := \hat{I}_{ij}(\mu)$.

1. The set of all the possible orders of $X$, $\pi(X)$, is the population $P$ of the sampling process (i.e. $P = \pi(X)$).
2. The matrix parameter under study is $\hat{I} = (\hat{I}_{ij})_{i,j \in X}$. For every pair of elements $i, j \in X$, $\hat{I}_{ij}$ is the estimated value of $I_{ij}$.
3. Each order $o \in \pi(X)$ is a sampling unit. Given a specific order $o = (o_1 \ldots o_n)$, the elements $o_r$ and $o_s$ are neighbors if $r = s + 1$ or $s = r + 1$. Under the assumption that $i$ and $j$ are consecutive or neighbors in the order $o$, the characteristics observed in each sampling unit are $\chi(o)_{ij} = \mu(Pred_j(o) \cup \{i\} \cup \{j\}) - \mu(Pred_j(o) \cup \{i\}) - \mu(Pred_j(o) \cup \{j\}) + \mu(Pred_j(o))$.
4. The estimated parameter $\hat{I}$ is the mean of the marginal contributions over the sample $M$, i.e. $\hat{I} = (\hat{I}_{ij})_{i,j \in X}$, where $\hat{I}_{ij} = \frac{1}{m_{ij}} \sum_{o \in M} \chi(o)_{ij}$, being $m_{ij}$ the amount of permutations among the $m$ possible sample orders in which $i$ and $j$ are neighbors. The total number of estimate all $I_{ij}$ is denoted as $Total Data = m(n-1)$.
5. The selection process considered to determine the sample $M$, takes any order $o \in \pi(X)$ with probability $\frac{1}{n!}$.

**Proposition 1.** *Considering the sample $M_{ij}$ with size $m_{ij}$, the estimator $\hat{I}_{ij}$ is unbiased, i.e. $E[\hat{I}_{ij}] = I_{ij}$. Moreover, its variance is given by*

---

**Algorithm 1** *ApproInteraction.*

1: **Input**: $(S, TotalData, \mu)$;
2: **Output**: $\hat{I} = \left(\hat{I}_{ij}\right)_{i,j \in S}$;
3: $m \leftarrow \frac{TotalData}{(n-1)}$;
4: $Cont \leftarrow 0$;
5: $\hat{I}_{ij} \leftarrow 0$; $m_{i,j} \leftarrow 0$ $\forall i, j \in S$;
6: **while** $(Cont < m)$ **do**
7:      $Cont \leftarrow Cont + 1$;
8:      Take $o = (o_{(1)}, \ldots, o_{(n)}) \in \pi(S)$ (with probability $\frac{1}{n!}$);
9:      **for** $(k = 1)$ **to** $(n - 1)$ **do**
10:          $m_{o_{(k)},o_{(k+1)}} \leftarrow m_{o_{(k)},o_{(k+1)}} + 1$;
11:          $m_{o_{(k+1)},o_{(k)}} \leftarrow m_{o_{(k)},o_{(k+1)}}$;
12:          $A \leftarrow \mu\left(Pred_{o_{(k)}}(o) \cup \{o_{(k)}\} \cup \{o_{(k+1)}\}\right)$;
13:          $B \leftarrow \mu\left(Pred_{o_{(k)}}(o) \cup \{o_{(k+1)}\}\right)$;
14:          $C \leftarrow \mu\left(Pred_{o_{(k)}}(o) \cup \{o_{(k)}\}\right)$;
15:          $D \leftarrow \mu\left(Pred_{o_{(k)}}(o)\right)$;
16:          $\chi(o)_{o_{(k)}o_{(k+1)}} \leftarrow A - B - C + D$;
17:          $\hat{I}_{o_{(k)}o_{(k+1)}} \leftarrow \hat{I}_{o_{(k)}o_{(k+1)}} + \chi(o)_{o_{(k)}o_{(k+1)}}$;
18:          $\hat{I}_{o_{(k+1)}o_{(k)}} \leftarrow \hat{I}_{o_{(k)}o_{(k+1)}}$;
19:      **end for**
20: **end while**
21: $\hat{I}_{ij} \leftarrow \frac{\hat{I}_{ij}}{m_{ij}}$, $\forall i, j \in S$;
22: **return** $(\hat{I})$;

---

$$Var[\hat{I}_{ij}] = \frac{\sigma^2}{m_{ij}}, \quad \text{where} \quad \sigma^2 = \frac{1}{(n-1)!} \sum_{o \in \pi(X \setminus \{i\})} \left(\chi(o)_{ij} - I_{ij}\right)^2$$

**Proof.** On the basis of random sampling processes [41], considering $\hat{I}_{ij}$ as a sample mean and $I_{ij}$ as a population mean, this proof is straightforward.

**Corollary 1.** *The estimator $\hat{I}_{ij}$ is consistent in probability, i.e.*

$$\lim_{m_{ij} \to \infty} p\left(|\hat{I}_{ij} - I_{ij}| > \epsilon\right) = 0, \quad \forall \epsilon > 0$$

**Proof.** Trivial.

**Corollary 2.** *Let $i, j \in X$. If $\forall K, K' \subset X \setminus \{i, j\}$ it holds that $\mu_{Kij} - \mu_{Ki} - \mu_{Kj} + \mu_K = \mu_{K'ij} - \mu_{K'i} - \mu_{K'j} + \mu_{K'}$, then $\hat{I}_{ij} = I_{ij}$ if $m_{ij} > 0$.*

**Proof.** Trivial.

Then we work in the determination of the sample size, $m_{ij}$. Given the parameters $e$ and $\alpha$, to fix $m_{ij}$, we have to guarantee that the error obtained in the estimation process is lower than $e$ with a probability greater than $(1 - \alpha)$. We apply the central limit theorem and the Proposition 1, so the property $\hat{I}_{ij} \sim N\left(I_{ij}, \frac{\sigma^2}{m_{ij}}\right)$ holds. Thus, if $m_{ij} \geq Z_{\frac{\alpha}{2}}^2 \frac{\sigma^2}{e^2}$, then $p\left(|\hat{I}_{ij} - I_{ij}| \leq e\right) \geq (1 - \alpha)$, where the value $Z_{\frac{\alpha}{2}}$ is such that $p\left(Z \geq Z_{\frac{\alpha}{2}}\right) = \frac{\alpha}{2}$ and $Z \sim N(0, 1)$.

On the other hand, as the value of $\sigma^2$ is unknown, it should be provided an upper bound of that value to determine the sample size $m_{ij}$. Then, we first define the minimum and the maximum values which can be reached by $\hat{I}_{ij}$: $\chi_{min} = \min_{\{o \in \pi(X \setminus \{i\}); i, j \in X\}} \chi(o)_{ij}$; and $\chi_{max} = \max_{\{o \in \pi(X \setminus \{i\}); i, j \in X\}} \chi(o)_{ij}$.

Let us remark that the maximum variance for any random variable bounded by two values is reached when that variable reaches one of these extreme values. This situation occurs with the same probability, $\frac{1}{2}$. In this case the two extreme values are $\chi_{min}$ and $\chi_{max}$, and the following inequality holds:

$$\sigma^2 \leq \frac{1}{2}\left(\chi_{max} - \frac{\chi_{max} + \chi_{min}}{2}\right)^2 + \frac{1}{2}\left(\chi_{min} - \frac{\chi_{max} + \chi_{min}}{2}\right)^2 = \frac{(\chi_{max} - \chi_{min})^2}{4}$$

Hence, the error is bounded by:

$$p\left(e \geq \sqrt{\frac{Z_{\frac{\alpha}{2}}^2}{m_{ij}} + \frac{(\chi_{max} - \chi_{min})^2}{4}}\right) \leq 1 - \alpha$$

## 4. Stratified estimation of the Interaction Index

There are some situations in which the *ApproInteraction* Algorithm may not be accurately enough, due to the large variance caused by simple random sampling. It occurs when there is a wide difference on the variance related to the different pairs of elements, depending on their position in the permutation considered. This source of variability can be reduced by the consideration of several strata. One of the most common techniques when considering Monte Carlo methods is the application of the stratification to estimate the expectations from a known population [24,25]. Then, following the idea introduced in [29], we propose a refinement of the algorithm introduced in Section 3 with the use of stratified random sampling, instead of simple random sampling. Hence, the performance of the original method is improved for cases in which there is a wide difference between the population variance and the variances on each stratum. Note that when these differences do not exist, it is not necessary to apply stratification. The key of the stratification is to divide the whole population into subpopulations, each of which is internally homogeneous and as heterogeneous as possible with respect to the other groups. To build these subgroups, we consider two of the most important sources of variation: the pairs of individuals, and the order in which each of them arrives to a coalition.

The distribution of the total sample used (named *TotalData* in *ApproInteraction* Algorithm) between different strata is performed with optimum allocation:

$$m_{ij\ell} = TotalData * \frac{\sigma_{ij\ell}^2}{\sum_{s,r \in X}\sum_{k=1}^{n-1} \sigma_{srk}^2}, \quad \forall i, j \in X, \ell = 1, \ldots, n-1. \tag{6}$$

To determine the optimum sample sizes per each stratum $P_{ij}^\ell$, the real value of the variance $\sigma_{ij\ell}^2$ has to be estimated. The proposed method is a two-stage approach, that divides the sample size $m_{ij\ell}$ into two parts, $m_{ij\ell} = m_{ij\ell}^{exp} + m_{ij\ell}^{st}$. Then, the different $m_{ij\ell}^{exp}$ samples are used to obtain unbiased estimations of the true variances $\sigma_{ij\ell}^2$. On the other hand, the $m_{ij\ell}^{st}$ are used to ensure that the final $m_{ij\ell}$ values are proportional to the resulting variance estimations.

The selection of $m_{ij\ell}^{exp}$ and $m_{ij\ell}^{st}$ has to be done depending on the problem. However, from the best of our knowledge, the effort needed to reach a good selection of these sizes is not worthy in general. Thus, we fix a balanced percentage of 50% per each part of the analysis for the computational results. Moreover, there are other features of the two-stage method that have to be set.

Next, we define a method to estimate $I_{ij}$ with a stratified process, named *StratifiedApproInteraction* (see Algorithm 2). In the sequel, the estimated value of $I_{ij}$ obtained with that stratified process will be denoted as $\hat{I}_{ij}^{st}$.

- **Stage 1.** *Pilot estimation*. For each position $\ell \in \{1, \ldots, n-1\}$ and $i, j \in X$:
  1. The sampling size in the experimental study is the same in all strata $P_{ij}^\ell$, i.e. $m_{ij\ell}^{exp} = \frac{TotalData}{n(n-1)^2}$.
  2. The characteristic observed in each sampling unit is $\chi(o)_{ij} = \mu\left(Pred_j(o) \cup \{i\} \cup \{j\}\right) - \mu\left(Pred_j(o) \cup \{i\}\right) - \mu\left(Pred_j(o) \cup \{j\}\right) + \mu\left(Pred_j(o)\right)$, where $i, j \in X$ and $o \in \pi(X\setminus\{i\})$.
  3. For every $i, j \in X$, $\ell = 1, \ldots, n-1$, we consider the stratum $P_{ij}^\ell$, which is the set of all the possible orders, $o \in \pi(X\setminus\{i\})$, in which the individual $j$ is in the $\ell$-th place.

**Algorithm 2** *StratifiedApproInteraction*.

1: **Input**: $(S, TotalData, \mu)$;
2: **Output**: $\hat{I}^{st} = \left(\hat{I}^{st}_{ij}\right)_{i,j \in S}$;
3: $\overline{T}^{exp}_{ij\ell} \leftarrow 0$; $\overline{T}^{st}_{ij\ell} \leftarrow 0$; $s^2_{ij\ell} \leftarrow 0, \forall i, j \in S, \forall \ell = 1, \ldots, n-1$;
4: **for** $(i = 1)$ **to** $(n)$ **do**
5:     **for** $(j = 1)$ **to** $(n)$ **do**
6:         **for** $(\ell = 1)$ **to** $(n-1)$ **do**
7:             $Cont1 \leftarrow 0$; $Cont2 \leftarrow 0$; $m^{exp}_{ij\ell} \leftarrow \frac{TotalData}{n(n-1)^2}$;
8:             **while** $\left(Cont1 < m^{exp}_{ij\ell}\right)$ **do**
9:                 $Cont1 \leftarrow Cont1 + 1$;
10:                Take $o = (o_{(1)}, \ldots, o_{(n-1)}) \in P^\ell_{ij}$ (with probability $\frac{1}{(n-2)!}$);
11:                $\chi(o)_{ij} \leftarrow \mu\left(Pred_j(o) \cup \{i\} \cup \{j\}\right) - \mu\left(Pred_j(o) \cup \{j\}\right) - \mu\left(Pred_j(o) \cup \{i\}\right) + \mu\left(Pred_j(o)\right)$;
12:                $\overline{T}^{exp}_{ij\ell} \leftarrow \overline{T}^{exp}_{ij\ell} + \chi(o)_{ij}$;
13:                $s^2_{ij\ell} \leftarrow s^2_{ij\ell} + (\chi(o)_{ij})^2$;
14:            **end while**
15:             $s^2_{ij\ell} \leftarrow \frac{1}{m^{exp}_{ij\ell}-1}(s^2_{ij\ell} - \frac{(\overline{T}^{exp}_{ij\ell})^2}{m^{exp}_{ij\ell}})$;
16:             $\overline{T}^{exp}_{ij\ell} \leftarrow \frac{\overline{T}^{exp}_{ij\ell}}{m^{exp}_{ij\ell}}$;
17:             $m^{st}_{ij\ell} \leftarrow \max\{0, \frac{TotalData * s^2_{ij\ell}}{\sum_{rt \in X} \sum_{k=1}^{n-1} s^2_{rtk}} - m^{exp}_{ij\ell}\}$;
18:             **while** $\left(Cont2 < m^{st}_{ij\ell}\right)$ **do**
19:                Take $o = (o_{(1)}, \ldots, o_{(n-1)}) \in P^\ell_{ij}$ (with probability $\frac{1}{(n-2)!}$);
20:                $\chi(o)_{ij} \leftarrow \mu\left(Pred_j(o) \cup \{i\} \cup \{j\}\right) - \mu\left(Pred_j(o) \cup \{j\}\right) - \mu\left(Pred_j(o) \cup \{i\}\right) + \mu\left(Pred_j(o)\right)$;
21:                $\overline{T}^{st}_{ij\ell} \leftarrow \overline{T}^{st}_{ij\ell} + \chi(o)_{ij}$;
22:            **end while**
23:             $\overline{T}^{st}_{ij\ell} \leftarrow \frac{\overline{T}^{st}_{ij\ell}}{m^{st}_{ij\ell}}$;
24:             $Cont2 \leftarrow Cont2 + 1$;
25:          **end for**
26:     **end for**
27: **end for**
28: $\hat{I}^{st}_{ij} \leftarrow \frac{1}{n-1} \sum_{\ell=1}^{n-1} \frac{(\overline{T}^{st}_{ij\ell} * m^{st}_{ij\ell}) + (\overline{T}^{exp}_{ij\ell} * m^{exp}_{ij\ell})}{(m^{st}_{ij\ell} + m^{exp}_{ij\ell})}, \forall i, j \in S$;
29: **return** $\left(\hat{I}^{st}\right)$;

    4. A process of random sampling with replacement is done in each stratum $P^\ell_{ij}$, obtaining a sample $M^{exp}_{ij\ell}$ whose size is $m^{exp}_{ij\ell}$. Let us note that the process done in each stratum to determine the corresponding sample takes any element of $P^\ell_{ij}$ with probability $\frac{1}{(n-2)!}$.

    5. For every $i, j \in X, \ell = 1, \ldots, n-1$, let $s^2_{ij\ell}$ be the unbiased sample variance of the marginal contributions over the sample $M^{exp}_{ij\ell}$, i.e. $s^2_{ij\ell} = \frac{1}{m^{exp}_{ij\ell}-1} \sum_{o \in M^{exp}_{ij\ell}} \left(\chi(o)_{ij} - \overline{T}^{exp}_{ij\ell}\right)^2$, where $\overline{T}^{exp}_{ij\ell}$ is the mean of $\chi(o)_{ij}$ over the sample $M^{exp}_{ij\ell}$.

- **Stage 2.** *Optimum allocation*:
    1. For each stratum, the total sampling size is proportional to the estimated variance $s^2_{ij\ell}$. Then, the sampling size $m^{st}_{ij\ell}$ is calculated for each stratum $P^\ell_{ij}$ as $m^{st}_{ij\ell} = m_{ij\ell} - m^{exp}_{ij\ell}$, where $m_{ij\ell} = TotalData * \frac{s^2_{ij\ell}}{\sum_{ij \in X} \sum_{k=1}^{n-1} s^2_{ijk}}$. Thus, in the second stage only the remaining $\left(m_{ij\ell} - m^{exp}_{ij\ell}\right)$ samples are randomly selected. If $m^{st}_{ij\ell} \leq 0$, no more samples are generated for the stratum $P^\ell_{ij}$, and all the other sample sizes $m_{ij\ell}$ are proportionally recalculated depending on $\left(m_{ij\ell} - m^{exp}_{ij\ell}\right)$.

2. The characteristic observed, the stratum and the process of random sampling are the same as in the Stage 1. The difference lies in the sampling size, $m_{ij\ell}^{st}$, and the sample, $M_{ij\ell}^{st}$.
3. Calculate $\overline{T}_{ij\ell}^{st}$ as the mean of $\chi(o)_{ij}$ over the sample $M_{ij\ell}^{st}$.
4. Calculate $\hat{I}_{ij}^{st} = \frac{(m_{ij\ell}^{st}*\overline{T}_{ij\ell}^{st})+(m_{ij\ell}^{exp}*\overline{T}_{ij\ell}^{exp})}{m_{ij\ell}^{exp}+m_{ij\ell}^{st}}$

**Proposition 2.** *The estimator $\hat{I}_{ij}^{st}$ obtained with StratifiedApproInteraction Algorithm is unbiased, i.e. $E[\hat{I}_{ij}^{st}] = I_{ij}$. Moreover, its variance is given by*

$$Var[\hat{I}_{ij}^{st}] = \frac{1}{(n-1)^2} \sum_{\ell=1}^{n-1} \frac{\sigma_{ij\ell}^2}{m_{ij\ell}}, \quad where \quad \sigma_{ij\ell}^2 = \frac{1}{(n-2)!} \sum_{o \in P_{ij}^\ell} \left(\chi(o)_{ij} - I_{ij}^\ell\right)^2 \quad (7)$$

**Proof.** On the basis of stratified random sampling processes [41], considering $\hat{I}_{ij}^{st}$ is a sample mean of a stratified sampling, and $I_{ij}$ is a population mean, this proof is straightforward.

**Corollary 3.** *The estimator $\hat{I}_{ij}^{st}$ is consistent in probability, i.e.*

$$\lim_{m_{ij} \to \infty} p\left(|\hat{I}_{ij}^{st} - I_{ij}| > \epsilon\right) = 0, \quad \forall \epsilon > 0$$

**Proof.** Trivial.

**Corollary 4.** *Let $i, j \in X$. If $\forall o, o' \in P_{ij}^\ell$ and $\forall \ell = 1, \ldots, n-1$ it holds that $\chi(o)_{ij} = \chi(o')_{ij}$, then $\hat{I}_{ij} = I_{ij}$ if $m_{ij} > 0$.*

**Proof.** Trivial.

To calculate the sample size $m_{ij\ell}$, the value $\sigma_{ij\ell}^2$ can be bounded by $\frac{(\chi_{max}-\chi_{min})^2}{4}$ (other bounds may be used instead). Then, the resulting theoretical error would be the same of that obtained with the method introduced in the previous section, if the sample sizes are similar.

## 5. Computational results

In this section we evaluate the performance of the new methods. To do so, we calculate the real value of $I_{ij}$ in several cases. Then, we compare the results provided by both algorithms *ApproInteraction* and *Stratified ApproInteraction* with those real values. Then, the first step is to calculate $I_{ij}$ in some well known examples. To do so, we consider the alternative calculation of the interaction index based on orders (see Theorem 1).

### 5.1. Some examples of calculating $I_{ij}$

Now we illustrate the calculation of the real value of $I_{ij}$ in several examples. We first define the values $A, B, C, D$ used in the different cases.

⋄ $A = \mu\left(Pred_j(o) \cup \{i\} \cup \{j\}\right)$
⋄ $B = \mu\left(Pred_j(o) \cup \{i\}\right)$
⋄ $C = \mu\left(Pred_j(o) \cup \{j\}\right)$
⋄ $D = \mu\left(Pred_j(o)\right)$

⋆ **Example 1. Symmetric stepped fuzzy measure**
We consider the function $\mu^1 : 2^X \longrightarrow [0, 1]$, defined in equation (8) (see Fig. 1), where $X = \{1, \ldots, 100\}$. It is a case of a symmetric fuzzy measure, in which all the individuals are equal, and the only important point is the size of $X$. In this situation, those subsets whose cardinality is relatively small, are almost worthless. Moreover, there is a value $t$ such that the obtained value is the same for every subset $S \subset X$ with $|S| \geq t$.

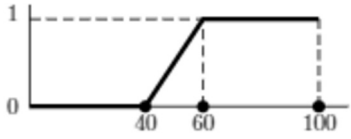

$$\mu^1(S) = \begin{cases} 0 & \text{if } 0 \leq s \leq 40 \\ 0.05(s-40) & \text{if } 40 \leq s \leq 60 \\ 1 & \text{otherwise} \end{cases}, \quad (8)$$

where $S \subseteq X$ and $s = |S|$.

Fig. 1. Fuzzy measure $\mu^1$.

Then we show the exact calculation of $I_{ij}$ regarding $\mu^1$. For every $i, j \in X$, $I_{ij}$ can be rewritten as

$$I_{ij}(\mu^1) = \sum_{o \in \pi(X \setminus \{i\})} \frac{1}{(n-1)!} \Big[ \mu^1\left(Pred_j(o) \cup \{i\} \cup \{j\}\right) - \mu^1\left(Pred_j(o) \cup \{i\}\right)$$

$$-\mu^1\left(Pred_j(o) \cup \{j\}\right) + \mu^1\left(Pred_j(o)\right) \Big] = \sum_{o \in \pi(X \setminus \{i\})} \frac{1}{(n-1)!} [A - B - C + D]$$

$$= \sum_{o \in \pi(X \setminus \{i\})} \frac{1}{(n-1)!} \left[ \chi(o)_{ij} \right]$$

- If $|Pred_j(o)| \leq 38$, then $A = B = C = D = 0$, so $\chi(o)_{ij} = 0$.
- If $|Pred_j(o)| = 39$, then $A = 0.05$ and $B = C = D = 0$, so $\chi(o)_{ij} = 0.05$.
- If $|Pred_j(o)| \in [40, 58]$, then $A = B + 0.05$ and $C = D + 0.05$. Then, $\chi(o)_{ij} = (B + 0.05) - B - (D + 0.05) + D = 0$.
- If $|Pred_j(o)| = 59$, then $A = B = C = 1$ and $D = 0.95$, so $\chi(o)_{ij} = 1 - 1 - 1 + 0.95 = -0.05$.
- If $|Pred_j(o)| \in [60, 98]$, then $A = B = C = D = 1$, so $\chi(o)_{ij} = 0$.

Then, the calculation of the interaction index $I_{ij}$ is:

$$I_{ij}(\mu^1) = \frac{1}{99!} \Bigg[ \sum_{\substack{o \in \pi(X \setminus \{i\})/ \\ |Pred_j(o)| \leq 38}} (0) + \sum_{\substack{o \in \pi(X \setminus \{i\})/ \\ |Pred_j(o)| = 39}} (0.05)$$

$$+ \sum_{\substack{o \in \pi(X \setminus \{i\})/ \\ |Pred_j(o)| \in [40,58]}} (0) + \sum_{\substack{o \in \pi(X \setminus \{i\}) \\ |Pred_j(o)| = 59}} (-0.05) + \sum_{\substack{o \in \pi(X \setminus \{i\}) \\ |Pred_j(o)| \in [60,98]}} (0) \Bigg]$$

On the other hand, the amount of orders in which the individual $j$ is in a particular $l$-th place is $\frac{(n-1)!}{n-1} = \frac{99!}{99}$, we can solve previous expression as

$$I_{ij}(\mu^1) = \frac{1}{99!} \left[ \frac{99!}{99} 0.05 + \frac{99!}{99} (-0.05) \right] = 0$$

Then, considering the fuzzy measure $\mu^1$ defined in equation (8), for every $i, j \in X$, the value of interaction index is $I_{ij}(\mu^1) = 0$.

★ **Example 2. Symmetric increasing fuzzy measure**

We consider the function $\mu^2 : 2^X \longrightarrow [0, 1]$, defined in equation (9) (see Fig. 2), where $X = \{1, \ldots, 100\}$. It is a case of a symmetric fuzzy measure, in which all the individuals are equal, and the only important thing is the size of $X$. The value assigned by $\mu^2$ to subsets with relatively small cardinal, grows slower than it does concerning bigger subsets.

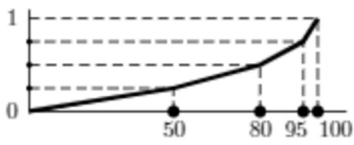

$$\mu^2(S) = \begin{cases} \frac{s}{200} & \text{if } 0 \leq s \leq 50 \\ 0.25 + \frac{s-50}{120} & \text{if } 50 \leq s \leq 80 \\ 0.5 + \frac{s-80}{60} & \text{if } 80 \leq s \leq 95 \\ 0.75 + \frac{s-95}{20} & \text{if } 95 \leq s \leq 100 \end{cases}, \quad (9)$$

where $S \subseteq X$ and $s = |S|$.

Fig. 2. Fuzzy measure $\mu^2$.

Then we show the exact calculation of $I_{ij}$ regarding $\mu^2$. For every $i, j \in X$, $I_{ij}$ can be rewritten as

$$I_{ij}\left(\mu^2\right) = \sum_{o \in \pi(X \setminus \{i\})} \frac{1}{(n-1)!} \left[\mu^2\left(Pred_j(o) \cup \{i\} \cup \{j\}\right) - \mu^2\left(Pred_j(o) \cup \{i\}\right)\right.$$
$$\left. - \mu^2\left(Pred_j(o) \cup \{j\}\right) + \mu^2\left(Pred_j(o)\right)\right] = \sum_{o \in \pi(X \setminus \{i\})} \frac{1}{(n-1)!} [A - B - C + D]$$
$$= \sum_{o \in \pi(X \setminus \{i\})} \frac{1}{(n-1)!} \left[\chi(o)_{ij}\right]$$

- If $|Pred_j(o)| \leq 48$, then $A = (0.005 + B)$ and $C = (0.005 + D)$, so $\chi(o)_{ij} = (B + 0.005) - B - (D + 0.005) + D = 0$.
- If $|Pred_j(o)| = 49$, then $A = \left(0.25 + \frac{1}{120}\right)$, $B = C = 0.25$ and $D = \left(0.25 - \frac{1}{200}\right)$, so $\chi(o)_{ij} = \left(0.25 + \frac{1}{120}\right) - 0.25 - \left(0.25 - \frac{1}{200}\right) + 0.25 = \frac{1}{120} - \frac{1}{200} = \frac{1}{300}$.
- If $|Pred_j(o)| \in [50, 78]$, then $A = \left(B + \frac{1}{120}\right)$ and $C = \left(D + \frac{1}{120}\right)$. Then, $\chi(o)_{ij} = \left(B + \frac{1}{120}\right) - B - \left(D + \frac{1}{120}\right) + D = 0$.
- If $|Pred_j(o)| = 79$, then $A = \left(0.5 + \frac{1}{60}\right)$, $B = C = 0.5$, and $D = \left(0.5 - \frac{1}{120}\right)$, so $\chi(o)_{ij} = \left(0.5 + \frac{1}{60}\right) - 0.5 - 0.5 + \left(0.5 - \frac{1}{120}\right) = \frac{1}{60} - \frac{1}{120} = \frac{1}{120}$.
- If $|Pred_j(o)| \in [80, 93]$, then $A = \left(B + \frac{1}{60}\right)$ and $C = \left(D + \frac{1}{60}\right)$, so $\chi(o)_{ij} = \left(B + \frac{1}{60}\right) - B - \left(D + \frac{1}{60}\right) + D = 0$.
- If $|Pred_j(o)| = 94$, then $A = \left(0.75 + \frac{1}{20}\right)$, $B = C = 0.75$, and $D = \left(0.75 - \frac{1}{60}\right)$, so $\chi(o)_{ij} = \left(0.75 + \frac{1}{20}\right) - 0.75 - 0.75 + \left(0.75 - \frac{1}{60}\right) = \frac{1}{20} - \frac{1}{60} = \frac{1}{30}$.
- If $|Pred_j(o)| \in [95, 98]$, then $A = \left(B + \frac{1}{20}\right)$ and $C = \left(D + \frac{1}{20}\right)$, so $\chi(o)_{ij} = \left(B + \frac{1}{20}\right) - B - \left(D + \frac{1}{20}\right) + D = 0$

Then, the calculation of the interaction index $I_{ij}$ is:

$$I_{ij}\left(\mu^2\right) = \frac{1}{99!}\left[\sum_{\substack{o \in \pi(X \setminus \{i\})/ \\ |Pred_j(o)| \leq 48}} (0) + \sum_{\substack{o \in \pi(X \setminus \{i\})/ \\ |Pred_j(o)|=49}} \left(\frac{1}{300}\right) + \sum_{\substack{o \in \pi(X \setminus \{i\})/ \\ |Pred_j(o)| \in [50,78]}} (0)\right.$$
$$\left. + \sum_{\substack{o \in \pi(X \setminus \{i\})/ \\ |Pred_j(o)|=79}} \left(\frac{1}{120}\right) + \sum_{\substack{o \in \pi(X \setminus \{i\})/ \\ |Pred_j(o)| \in [80,93]}} (0) + \sum_{\substack{o \in \pi(X \setminus \{i\})/ \\ |Pred_j(o)|=94}} \left(\frac{1}{30}\right) + \sum_{\substack{o \in \pi(X \setminus \{i\})/ \\ |Pred_j(o)| \in [95,100]}} (0)\right]$$

On the other hand, the amount of orders in which the individual $j$ is in a particular $l$-th place is $\frac{(n-1)!}{n-1} = \frac{99!}{99}$, we can solve previous expression as

$$I_{ij}\left(\mu^2\right) = \frac{1}{99!}\left[\frac{99!}{99}\frac{1}{300} + \frac{99!}{99}\frac{1}{120} + \frac{99!}{99}\frac{1}{30}\right] = \frac{1}{99}\left(\frac{1}{300} + \frac{1}{120} + \frac{1}{30}\right) = \frac{1}{2200}.$$

Then, considering the fuzzy measure $\mu^2$ defined in equation (9), for every $i, j \in X$, the value of the interaction index is $I_{ij}\left(\mu^2\right) = \frac{1}{2200}$.

* **Example 3. Fuzzy measure with maximum weight**
We consider the function $\mu^3 : 2^X \longrightarrow [0,1]$, defined in equation (10), where $X = \{1, \ldots, 100\}$ and $S \subseteq X$. This function is a fuzzy measure in which all the individuals have an assigned weight. In this situation, the result related to a particular subset, is the maximum value of the weights assigned to its individuals. It is known as the airport problem.

$$\mu^3(S) = \max_{i \in S}\{w_i\}, \text{ where } \begin{cases} w_i = 0.1, & \forall i \in \{1, \ldots, 8\} \\ w_i = 0.2, & \forall i \in \{9, \ldots, 20\} \\ w_i = 0.3, & \forall i \in \{21, \ldots, 26\} \\ w_i = 0.4, & \forall i \in \{27, \ldots, 40\} \\ w_i = 0.5, & \forall i \in \{41, \ldots, 48\} \\ w_i = 0.6, & \forall i \in \{49, \ldots, 57\} \\ w_i = 0.7, & \forall i \in \{58, \ldots, 70\} \\ w_i = 0.8, & \forall i \in \{71, \ldots, 80\} \\ w_i = 0.9, & \forall i \in \{81, \ldots, 90\} \\ w_i = 1, & \forall i \in \{91, \ldots, 100\} \end{cases} \qquad (10)$$

Then we show the exact calculation of $I_{ij}$ regarding $\mu^3$. For every $i, j \in X$, $I_{ij}$ can be rewritten as

$$I_{ij}(\mu^3) = \sum_{o \in \pi(X \setminus \{i\})} \frac{1}{(n-1)!} \left[ \mu^3\left(Pred_j(o) \cup \{i\} \cup \{j\}\right) - \mu^3\left(Pred_j(o) \cup \{i\}\right) \right.$$
$$\left. -\mu^3\left(Pred_j(o) \cup \{j\}\right) + \mu^3\left(Pred_j(o)\right) \right] = \sum_{o \in \pi(X \setminus \{i\})} \frac{1}{(n-1)!} [A - B - C + D]$$
$$= \sum_{o \in \pi(X \setminus \{i\})} \frac{1}{(n-1)!} \left[ \chi(o)_{ij} \right]$$

We define $Min_{ij} = \min\{w_i, w_j\}$; and, without loss of generality, we assume $w_i \geq w_j$. Then, $I_{ij}(\mu^3) = \sum_{o \in \pi(X \setminus \{i\})} \frac{1}{(n-1)!} [A - B - C + D] = -\sum_{o \in \pi(X \setminus \{i\})} \frac{1}{(n-1)!} [C - D]$.

The only difference between $A$ and $B$ is the element $\{j\}$, so both values are equal. Then, under the assumption that $w_i \geq w_j$, we have that $I_{ij}(\mu^3) = -\sum_{o \in \pi(X \setminus \{i\})} \frac{1}{(n-1)!} [C - D]$ is minus the Shapley index of the individual $j$ on $\mu^3$, when $i$ 'does not exist'. Concretely, $I_{ij}(\mu^3) = -Sh_j(\mu^{3i})$, where $\mu^{3i}(S) = \mu^3(S), \forall S \subseteq X^i = X \setminus \{i\}$. The calculation of the Shapley index of these phased fuzzy measures, $\mu^3$ and $\mu^{3i}$, has an explicit calculation, as it can be seen in [28]. Note that, for every pair of individuals $i, j$, the value $Min_{ij}$ only depends on the weights $w_i$ and $w_j$, so there is no need to calculate the value of $I_{ij}$ for every feasible pair, as the complete casuistry is reduced to the lines below.

– If $Min_{ij} = 0.1$, then $I_{ij}(\mu^3) = -\left[\frac{0.1}{99}\right]$
– If $Min_{ij} = 0.2$, then $I_{ij}(\mu^3) = -\left[\frac{0.1}{99} + \frac{0.1}{91}\right]$
– If $Min_{ij} = 0.3$, then $I_{ij}(\mu^3) = -\left[\frac{0.1}{99} + \frac{0.1}{91} + \frac{0.1}{79}\right]$
– If $Min_{ij} = 0.4$, then $I_{ij}(\mu^3) = -\left[\frac{0.1}{99} + \frac{0.1}{91} + \frac{0.1}{79} + \frac{0.1}{73}\right]$
– If $Min_{ij} = 0.5$, then $I_{ij}(\mu^3) = -\left[\frac{0.1}{99} + \frac{0.1}{91} + \frac{0.1}{79} + \frac{0.1}{73} + \frac{0.1}{53}\right]$
– If $Min_{ij} = 0.6$, then $I_{ij}(\mu^3) = -\left[\frac{0.1}{99} + \frac{0.1}{91} + \frac{0.1}{79} + \frac{0.1}{73} + \frac{0.1}{53} + \frac{0.1}{51}\right]$
– If $Min_{ij} = 0.7$, then $I_{ij}(\mu^3) = -\left[\frac{0.1}{99} + \frac{0.1}{91} + \frac{0.1}{79} + \frac{0.1}{73} + \frac{0.1}{53} + \frac{0.1}{51} + \frac{0.1}{42}\right]$
– If $Min_{ij} = 0.8$, then $I_{ij}(\mu^3) = -\left[\frac{0.1}{99} + \frac{0.1}{91} + \frac{0.1}{79} + \frac{0.1}{73} + \frac{0.1}{53} + \frac{0.1}{51} + \frac{0.1}{42} + \frac{1}{29}\right]$
– If $Min_{ij} = 0.9$, then $I_{ij}(\mu^3) = -\left[\frac{0.1}{99} + \frac{0.1}{91} + \frac{0.1}{79} + \frac{0.1}{73} + \frac{0.1}{53} + \frac{0.1}{51} + \frac{0.1}{42} + \frac{1}{29} + \frac{1}{19}\right]$
– If $Min_{ij} = 1$, then $I_{ij}(\mu^3) = -\left[\frac{0.1}{99} + \frac{0.1}{91} + \frac{0.1}{79} + \frac{0.1}{73} + \frac{0.1}{53} + \frac{0.1}{51} + \frac{0.1}{42} + \frac{1}{29} + \frac{1}{19} + \frac{1}{9}\right]$

For example, if $i = 4, j = 35$, $Min_{ij} = \min\{w_4, w_{35}\} = \min\{0.1, 0.4\} = 0.1$, which corresponds to the first case, so $I_{4,35}(\mu^3) = -\left[\frac{0.1}{99}\right]$.

⋆ **Example 4. Bi-cluster fuzzy measure**

We consider the function $\mu^4 : 2^X \to [0, 1]$, defined in equation (11) with $X = \{1, \ldots, 100\}$ and $S \subseteq X$. It is an example of fuzzy measure in which the elements are organized into two groups, and the value related to a set is the amount of pairs (each one belonging to one group) which could be defined following any rule, the popular shoes problem.

$$\mu^4(S) = 0.02 \min\{n_1(S), n_2(S)\}, \quad \text{where } n_i(S) = |C_S^i|, \ i = 1, 2 \tag{11}$$

$$C_S^1 = S \cap \{1, \ldots, 50\} \text{ and } C_S^2 = S \cap \{51, \ldots, 100\}$$

Then we show the exact calculation of $I_{ij}$ regarding $\mu^4$. For every $i, j \in X$, $I_{ij}$ can be rewritten as

$$I_{ij}\left(\mu^4\right) = \sum_{o \in \pi(X \setminus \{i\})} \frac{1}{(n-1)!} [\mu^4\left(Pred_j(o) \cup \{i\} \cup \{j\}\right) - \mu^4\left(Pred_j(o) \cup \{i\}\right)$$

$$-\mu^4\left(Pred_j(o) \cup \{j\}\right) + \mu^4\left(Pred_j(o)\right)] = \sum_{o \in \pi(X \setminus \{i\})} \frac{1}{(n-1)!} [A - B - C + D]$$

$$= \sum_{o \in \pi(X \setminus \{i\})} \frac{1}{(n-1)!} \left[\chi(o)_{ij}\right]$$

Then, we distinguish several cases:

1. If $i, j \in \{1, \ldots, 50\} = C_X^1$
   (a) If $n_1\left(Pred_j(o)\right) \geq n_2\left(Pred_j(o)\right)$, then, $A = B = C = D$, so $\chi(o)_{ij} = 0$.
   (b) If $n_1\left(Pred_j(o)\right) \leq \left(n_2\left(Pred_j(o)\right) - 2\right)$, then, $A = (B + 0.02)$ and $C = (D + 0.02)$, so $\chi(o)_{ij} = 0$.
   (c) If $n_1\left(Pred_j(o)\right) = \left(n_2\left(Pred_j(o)\right) - 1\right)$, then, $A = B = C = (D + 0.02)$, so $\chi(o)_{ij} = -0.02$.
   
   In conclusion, if $\exists c \in \mathbb{N}$ such that $|Pred_j(o)| = 2c$, then $I_{ij}\left(\mu^4\right) = 0$. Otherwise, if $\exists c \in \mathbb{N}$ such that $|Pred_j(o)| = 2c - 1$, then the probability that $I_{ij}\left(\mu^4\right) = -0.02$ is the same probability of obtaining $c$ elements of type $C_S^2$ and $(c-1)$ elements of type $C_S^1$, when there are 50 elements of type 2, and 48 elements of type 1. This probability is calculated by means of the hypergeometric distribution, where $\mathcal{P}_\mathcal{H}(a_1, a_2, A_1, A_2) = \frac{\binom{A_1}{a_1}\binom{A_2}{a_2}}{\binom{A_1+A_2}{a_1+a_2}}$. Furthermore, there are $\frac{(n-1)!}{n-1} = \frac{99!}{99}$ orders $o$ for which the property $|Pred_j(o)| = k$ holds. Then, $I_{ij}\left(\mu^4\right) = \sum_{o \in \pi(S \setminus \{i\})} \frac{1}{99!} \left[\chi(o)_{ij}\right] = \frac{1}{99!} \sum_{k=0}^{98} \sum_{\substack{o \in \pi(S \setminus \{i\})/ \\ |Pred_j(o)|=k}} \left[\chi(o)_{ij}\right] =$
   
   $\sum_{c=1}^{49} \frac{99!}{99} \mathcal{P}_\mathcal{H}(c, c-1, 50, 48) (-0.02) = \frac{-0.02}{99} \sum_{c=1}^{49} \frac{\binom{50}{c}\binom{48}{c-1}}{\binom{98}{2c-1}} = -0.00215602304$.

2. If $i, j \in \{51, \ldots, 100\}$, the process is analogous to point 1., interchanging the role of the classes $C_S^1$ and $C_S^2$. The result is the same.

3. Otherwise, $i \in C_X^1$ and $j \in C_X^2$, or vice versa. Without loss of generality, we assume $i \in \{1, \ldots, 50\} = C_X^1$, and $j \in \{51, \ldots, 100\} = C_X^2$.
   (a) If $n_1\left(Pred_j(o)\right) > n_2\left(Pred_j(o)\right)$, then $A = C$ and $B = D$, so $\chi(o)_{ij} = 0$.
   (b) If $n_1\left(Pred_j(o)\right) < n_2\left(Pred_j(o)\right)$, then $A = B$ and $C = D$, so $\chi(o)_{ij} = 0$.
   (c) If $n_1\left(Pred_j(o)\right) = n_2\left(Pred_j(o)\right)$, then $A = (B + 0.02) = (C + 0.02) = (D + 0.02)$, so $\chi(o)_{ij} = 0.02$.
   
   In conclusion, if $\exists c \in \mathbb{N}$ such that $|Pred_j(o)| = 2c - 1$, then $I_{ij}\left(\mu^4\right) = 0$. Otherwise, if $\exists c \in \mathbb{N}$ such that $|Pred_j(o)| = 2c$, then the probability that $I_{ij}\left(\mu^4\right) = 0.02$ is the same probability of obtaining $c$ elements of type $C_S^1$ and $c$ elements of type $C_S^2$, when there are 49 elements of type 2, and 49 elements of type 1. This probability is calculated by means of the hypergeometric distribution, where $\mathcal{P}_\mathcal{H}(a_1, a_2, A_1, A_2) = \frac{\binom{A_1}{a_1}\binom{A_2}{a_2}}{\binom{A_1+A_2}{a_1+a_2}}$. Furthermore, there are $\frac{(n-1)!}{n-1} = \frac{99!}{99}$ orders $o$ for which the property $|Pred_j(o)| = k$ holds. Then, $I_{ij}\left(\mu^4\right) = \sum_{o \in \pi(X \setminus \{i\})} \frac{1}{99!} \left[\chi(o)_{ij}\right] = \frac{1}{99!} \sum_{k=0}^{98} \sum_{\substack{o \in \pi(X \setminus \{i\})/ \\ |Pred_j(o)|=k}} \left[\chi(o)_{ij}\right] =$
   
   $\sum_{c=1}^{49} \frac{99!}{99} \mathcal{P}_\mathcal{H}(c, c, 49, 49) (0.02) = \frac{0.02}{99} \sum_{c=1}^{49} \frac{\binom{49}{c}\binom{49}{c}}{\binom{98}{2c}} = 0.00251290258037098$.

Table 1
Example 1. Symmetric stepped fuzzy measure. Errors made by the algorithms *ApproInteraction* and *StratifiedApproInteraction*.

| $AvgData$ | $e_{max}$ | $e_{av}$ | $e_{max}^{st}$ | $e_{av}^{st}$ | $e_{th}$ |
|---|---|---|---|---|---|
| $10^3$ | 8.08081 | 1.81791 | 0 | 0 | 5.20195 |
| $10^4$ | 2.74474 | 0.57115 | 0 | 0 | 1.64500 |
| $10^5$ | 0.81294 | 0.17876 | 0 | 0 | 0.52019 |
| $10^6$ | 0.26205 | 0.05690 | 0 | 0 | 0.16450 |
| $10^7$ | 0.08922 | 0.01815 | 0 | 0 | 0.05202 |

Table 2
Example 2. Symmetric increasing fuzzy measure. Errors made by the algorithms *ApproInteraction* and *StratifiedApproInteraction*.

| $AvgData$ | $e_{max}$ | $e_{av}$ | $e_{max}^{st}$ | $e_{av}^{st}$ | $e_{th}$ |
|---|---|---|---|---|---|
| $10^3$ | 4.10515 | 0.85948 | 0 | 0 | 1.73398 |
| $10^4$ | 1.28224 | 0.27403 | 0 | 0 | 0.54833 |
| $10^5$ | 0.39009 | 0.08637 | 0 | 0 | 0.17340 |
| $10^6$ | 0.13403 | 0.02755 | 0 | 0 | 0.05483 |
| $10^7$ | 0.03908 | 0.00887 | 0 | 0 | 0.01734 |

### 5.2. Testing the performance of the algorithms

Finally we evaluate our algorithms. To do so, we consider the real values of $I_{ij}$ calculated in polynomial time for fuzzy measures $\mu^1$, $\mu^2$, $\mu^3$ and $\mu^4$ in previous section. Then, we estimate the corresponding value of the interaction index with the new proposed methods. These examples are used to quantify the error of the sampling processes proposed in *ApproInteraction* Algorithm and *StratifiedApproInteraction* Algorithm. Let us note that, in the absence of other methods to approximate the value of the interaction index, we measure the goodness of the algorithms here proposed by comparing the obtained result with the real value of those examples. We show that both methods are so accurate, so the obtained results are exact and reliable. Hence, in Tables 1-4 we show the precision of the algorithms *ApproInteraction* and *StratifiedApproInteraction* for the examples previously introduced. To be clearer, we consider the scale 1 : 0.001, it is an error of 0.001 will be represented with the value 1. In these tables we show the average sample size for each $I_{ij}$, ($AvgData = TotalData/\frac{n(n-1)}{2}$). Then, $e_{max}$ and $e_{av}$ are, respectively, the maximum and average error considering all the estimations of $I_{ij}$ obtained with the *ApproInteraction* Algorithm. Similarly, $e_{max}^{st}$ and $e_{av}^{st}$ are, respectively, the maximum and average error considering all the estimations of $I_{ij}$ obtained with the *StratifiedApproInteraction* Algorithm, and $e_{th}$ is the bound of the theoretical error when considering $\alpha = 0.001$, $p\left(|I_{ij} - \hat{I}_{ij}|\right) \geq (1 - \alpha) = 0.999$. The calculation of the theoretical error depends on the case addressed:

- Example 1, $\mu^1$: $\chi_{max}^j = 0.05$ and $\chi_{min}^j = -0.05$.
- Example 2, $\mu^2$: $\chi_{max}^j = \frac{1}{30}$ and $\chi_{min}^j = 0$.
- Example 3, $\mu^3$: $\chi_{max}^j = 0$ and $\chi_{min}^j = -1$.
- Example 4, $\mu^4$: $\chi_{max}^j = 0.02$ and $\chi_{min}^j = -0.02$.

Note that, for those fuzzy measures in which the value $\chi(o)_{ij}$ varies a lot and is largely dependent on the values $i, j \in X$ and the position $\ell$ that these values have in the order $o$ considered, the *StratifiedApproInteraction* Algorithm will perform better than *ApproInteraction* Algorithm. This situation occurs in Examples 1, 2, 3. Otherwise, if the value $\chi(o)_{ij}$ does not vary a lot, then both methods perform just as well, as it can be seen in the Example 4.

Note that the *StratifiedApproInteraction* Algorithm does not make any error for Example 1 and Example 2, because of Corollary 4. Below we show the analysis of each example in the corresponding table. For a better visualization of the data below, we consider the scale 1 : 0.0001, i.e., the values showed in the following tables correspond with the obtained result multiplied by 1000.

Table 3
Example 3. Fuzzy measure with maximum weight. Errors made by the algorithms *ApproInteraction* and *StratifiedApproInteraction*.

| $AvgData$ | $e_{max}$ | $e_{av}$ | $e_{max}^{st}$ | $e_{av}^{st}$ | $e_{th}$ |
|---|---|---|---|---|---|
| $10^3$ | 115.92435 | 11.16639 | 45.47629 | 3.31031 | 52.01947 |
| $10^4$ | 43.03014 | 3.46081 | 4.29888 | 0.33309 | 16.45000 |
| $10^5$ | 10.06615 | 1.08402 | 1.03843 | 0.08281 | 5.20195 |
| $10^6$ | 3.69374 | 0.33845 | 0.26446 | 0.02507 | 1.64500 |
| $10^7$ | 0.93719 | 0.10004 | 0.09780 | 0.00795 | 0.52019 |

Table 4
Example 4. Bi-cluster fuzzy measure. Errors made by the algorithms *ApproInteraction* and *StratifiedApproInteraction*.

| $AvgData$ | $e_{max}$ | $e_{av}$ | $e_{max}^{st}$ | $e_{av}^{st}$ | $e_{th}$ |
|---|---|---|---|---|---|
| $10^3$ | 8.41007 | 1.61541 | 6.74519 | 1.40705 | 2.08078 |
| $10^4$ | 2.71840 | 0.51198 | 1.90963 | 0.42412 | 0.65800 |
| $10^5$ | 0.73839 | 0.16361 | 0.63031 | 0.11466 | 0.20808 |
| $10^6$ | 0.28346 | 0.05153 | 0.16600 | 0.03596 | 0.06580 |
| $10^7$ | 0.08438 | 0.01921 | 0.05536 | 0.01183 | 0.02081 |

## 6. Discussion and conclusions

Fuzzy measures comprise an important field with lot of applications in many disciplines, as modeling interactions between individuals. Handling fuzzy measures is not easy, nor does giving an interpretation to them. One of the most important contributions in this field is the Shapley value [21], an index defined to quantify the importance of each singleton in a coalition.

As a generalization of it, Murofushi and Soneda proposed an interaction index [36] able to model the interaction between pairs of elements. This index, an essential contribution on the interpretation of fuzzy measures and with base on some elements borrowed from multi-attribute utility theory, can be understood as the average value obtained when considering two elements together, taking into account all the coalitions. Then, Grabisch proposed an extension of the interaction index providing a characterization to deal with the relations inherent in the elements of any set [19]. It is the cornerstone of this paper.

To face the complex calculation of the Shapley index, in [28] it was proposed a sampling method to estimate this value in polynomial time. Shortly after, in [29], it was proposed a refinement of this method, based on a process of stratified random sampling with minimum allocation (highly recommendable to cope with situations which are prone to exhibit high variability on the individual or on the order considered).

Following the philosophy of those papers, we propose two sampling methods to estimate the interaction index. To carry on with it, inspired by the characterization of the Shapley value based on orders, we provide an alternative definition of the interaction index based in orders. This characterization is used to demonstrate that, in some cases, the interaction index can be calculated in polynomial time. Then, we provide two algorithms based on sampling to estimate it. Both methods have desirable properties; specifically, the obtained value is unbiased. The method proposed in Section 3 is based on random sampling. This algorithm is quite simple and it is easy to compute; it is a good pre-step to understand and process the stratified version described in Section 4.

To test the goodness of our methods we provide some computational results. We compare the results obtained with *ApproInteraction* Algorithm and *StratifiedApproInteraction* Algorithm with the real value of some cases of fuzzy measures for which we calculate the real value of the corresponding interaction index 5.1. As it can be seen in Section 5, both algorithms perform very well.

We would like to stress that, although we focused on the interaction index $I_{ij}$, the stratified procedure can be generalized in a quite similar way to estimate other values. Particularly, the estimation of $I_T(\mu)$ can be done just modifying the Algorithm 2.

**Declaration of competing interest**

The authors declare the following financial interests/personal relationships which may be considered as potential competing interests: Daniel Gómez reports financial support was provided by Government of Spain. Inmaculada Gutiérrez reports financial support was provided by Complutense University of Madrid.

**Data availability**

No data was used for the research described in the article.

**Acknowledgements**

Funding: this research has been partially supported by the Government of Spain, Grant Plan Nacional de I+D+i, [PID2020-116884GB-I00, PR108/20-28, PGC2018096509-B-I00] and the Complutense University of Madrid, [CT17/17 - CT18/17].